\documentclass[12pt, reqno]{elsart}
\usepackage{amssymb}
\usepackage{amsmath}  	% for multi-line subscripts etc.
\usepackage{xspace}		% this makes text macros not annoying to use
\usepackage{amscd}		% for drawing simple commutative diagrams
\usepackage{natbib}
\usepackage{bm}

\newtheorem{theorem}{Theorem}[section]
\newtheorem{lemma}[theorem]{Lemma}
\newtheorem{corollary}[theorem]{Corollary}
\newtheorem{proposition}[theorem]{Proposition}
\newtheorem{question}[theorem]{Question}

 \theoremstyle{definition}

\newtheorem{definition}[theorem]{Definition}

\newtheorem{remark}[theorem]{Remark}

 \theoremstyle{remark}

\newcommand{\C}{\mathbb{C}}
\newcommand{\R}{\mathbb{R}}

\newcommand{\N}{\mathbb{N}}

\newcommand{\mf}[1]{\mathfrak{#1}}

\newcommand{\ra}{\rightarrow}

\newcommand{\trdeg}{\operatorname{trdeg}}
\newcommand{\GK}{\operatorname{GKdim}}

\newcommand{\Krull}{\operatorname{Kdim}}
\newcommand{\Hom}{\operatorname{Hom}}

\newcommand{\End}{\operatorname{End}}

\newcommand{\spec}{\operatorname{spec}}
\newcommand{\mspec}{\operatorname{mspec}}
\newcommand{\ideal}{\lhd}

\newcommand{\set}[1]{\{  #1  \} }
\newcommand{\setmid}[2]{\{\,  #1 \mid #2 \,\} }

\newcommand{\Jac}{J}

\newcommand{\extcent}[1]{C(#1)}
\newcommand{\centclose}[1]{#1 \extcent{#1}}
\newcommand{\mart}[1]{Q_{r}(#1)}

\newcommand{\tcap}{\textstyle\bigcap\limits}

\newcommand{\fg}{finitely generated\xspace}
\newcommand{\PI}{\textup{PI}\xspace}
\newcommand{\ji}{just infinite\xspace}
\newcommand{\JI}{Just Infinite\xspace}
\newcommand{\sji}{stably just infinite\xspace}

\begin{document}
\begin{frontmatter}
\title{Stably Just Infinite Rings}

\author{Jason Bell\thanksref{JB}}
\thanks[JB]{The first author thanks NSERC for its generous support.}
\address{
	Department of Mathematics \\
	Simon Fraser University \\
	Burnaby, BC V5A 1S6. Canada	}
\ead{jpb@math.sfu.ca}

\author{John Farina\corauthref{JF}}
\corauth[JF]{Corresponding author.}
\address{
	Department of Mathematics \\
	University of California, San Diego \\
	La Jolla, CA  92093 }
\ead{jfarina@math.ucsd.edu}

\author{Cayley Pendergrass-Rice\thanksref{CP}}
\thanks[CP]{The third author received partial support from a grant from the Hewlett Mellon Fund for Faculty Development at Albion College.}
\address{
	Department of Mathematics \\
	Albion College \\
	Albion, MI 49224 }
\ead{cpendergrass@albion.edu}

\date{}

\begin{abstract}
We study \ji algebras which remain so upon extension of scalars by arbitrary field extensions.  Such rings are called \emph{\sji.}  We show that \ji rings over algebraically closed fields are \sji provided that the ring is either right noetherian (\ref{cor:noetherian-non-pi-ji-implies-sji}) or countably generated over a large field (\ref{cor:cg-non-PI-ji-over-big-fields-are-sji}).  We give examples to show that, over countable fields, a \ji algebra which is either affine or non-noetherian need not remain \ji under extension of scalars.  We also give a concrete classification of \PI \sji rings (\ref{thm:pi-sji-characterization}) and give two characterizations of non-\PI \sji rings in terms of Martindale's \emph{extended center} (\ref{thm:non-pi-sji-iff-centrally-closed}, \ref{thm:sji-iff-tensor-with-extcent}). 
\end{abstract}


\end{frontmatter}

\section{Introduction}

Simple rings are so-called because their two-sided ideal structure is as simple as possible.  A natural way of generalizing the notion of a simple ring is to allow the presence of two-sided ideals, but to insist that the nonzero ones be very ``large.''  Specifically, we consider algebras for which all nonzero ideals have finite codimension.  We begin, as always, with a definition.         
\begin{definition}\em{
A $k$-algebra $A$ is called \emph{\ji dimensional}, or \emph{\ji} for short, if $\dim_{k}(A) = \infty$ and each of its nonzero two-sided ideals has finite codimension.}
\end{definition}
The nomenclature which has been adopted is to call such rings \emph{projectively simple} in case the ring is $\N$-graded.   There are many examples of just infinite algebras.  Aside from any infinite dimensional simple ring, the most immediate example is a polynomial ring in one variable $k[x]$.  Reichstein, Rogalski, and Zhang \cite{ReichsteinRogalskiZhang:Projectively-simple-rings} give non-PI examples constructed using abelian varieties and Bartholdi \cite{Bartholdi:Branch-Rings-Thinned} constructs interesting \ji rings using groups acting on infinite trees.  

Much is known for just infinite algebras over uncountable fields; for example, Farkas and Small \cite{FarkasSmall:Algebras-which-are-nearly} show they are either primitive, have nonzero Jacobson radical, or satisfy a polynomial identity.  Over countable fields much less is known and many of the basic questions that are easily answered over uncountable fields are still open.  For this reason we are interested in when we can extend scalars in a just infinite algebra and remain just infinite.  We thus give the following definition. 

\begin{definition}\em{
A \ji $k$-algebra $A$ is called \emph{\sji} if $A\otimes_{k}K$ is \ji over $K$, for every field extension $K/k$.}
\end{definition} 

Just infinite algebras, even affine ones, may behave rather poorly upon extension of scalars.  For example, $\C[x]$ is an affine \ji $\R$-algebra, but $\C[x]\otimes_{\R}\C$ contains $\C\otimes_{\R}\C$, which is not a domain.  Thus $\C[x]\otimes_{\R}\C$ is not even prime and hence $\C[x]$ fails to be \ji over $\R$ \cite{FarinaPendergrass:A-Few-Properties-of-Just-Infinite}.  This example shows that in case the ground field is not algebraically closed, things can go horribly awry even upon extending scalars by a finite dimensional field extension!   Bergman proved the tensor product of any two prime $k$-algebras is again prime when $k$ is algebraically closed \cite[Proposition 17.2]{Bergman:Zero-divisors-in-tensor-products}.  It is therefore reasonable to restrict our attention to the case that the ground field is algebraically closed.  In this way we can, at the very least, avoid the unpleasantness of having a \ji algebra $A$ for which extension of scalars results in a ring which fails to be prime.

We summarize our main results:  We show that right noetherian \ji algebras are well-behaved under extension of scalars.  

\begin{theorem}\label{thm:non-pi-noetherian-ji-implies-sji} Let $k$ be an algebraically closed field and let $A$ be a non-PI right noetherian \ji $k$-algebra.  Then $A$ is stably \ji.
\end{theorem}

In addition, we show that the hypothesis that $A$ is a right noetherian ring cannot be weakened when the ground field is countable---even if the algebra is \fg.

\begin{theorem}\label{thm:affine-ji-non-sji-rings-over-countable-fields-exist} Let $k$ be a countable algebraically closed field.  Then there exists an affine non-simple right Goldie \ji $k$-algebra which is not \sji.
\end{theorem}

It is necessary to take $k$ to be a countable algebraically closed field in this theorem.  In fact, we show that if $A$ is a countably generated  \ji $k$-algebra over an uncountable algebraically closed field, then $A$ is \sji.  

\begin{theorem}\label{thm:cg-ji-over-big-fields-are-sji}
Let $k$ be an uncountable algebraically closed field and let $A$ be a countably generated \ji $k$-algebra.  Then $A$ is \sji.	
\end{theorem}

We remark that this gives an analog of a result of Reichstein et al. \cite[Lemma 1.8]{ReichsteinRogalskiZhang:Projectively-simple-rings} for ungraded rings.  To complete our study of \sji algebras, we look at the \PI case and show that a \PI \ji $k$-algebra $A$ is \sji if and only if the center of the quotient division algebra of $A$ is a \fg field extension of $k$ of transcendence degree $1$.

This paper is organized as follows.  In \S \ref{sec:ji-rings} and \S \ref{sec:central-closure}, we give some useful facts about \ji algebras and extended centers, which we later use in our consideration of \sji algebras.  In \S \ref{sec:noetherian-ji-rings}, we prove Theorem \ref{thm:non-pi-noetherian-ji-implies-sji}.  In \S \ref{sec:pi-ji-rings}, we present our characterization of \PI \sji algebras.  In \S 
\ref{sec:uncountable-fields}, we turn our attention to \ji algebras over uncountable fields and prove Theorem \ref{thm:cg-ji-over-big-fields-are-sji}.  Finally, in \S \ref{sec:examples}, we construct interesting \ji algebras that are not \sji to obtain Theorem \ref{thm:affine-ji-non-sji-rings-over-countable-fields-exist}.

All rings under consideration will be associative algebras, with $1$, over a field $k$, and homomorphisms are assumed to be unital.

\section{Just Infinite Algebras}\label{sec:ji-rings}

Reichstein et al. \cite{ReichsteinRogalskiZhang:Projectively-simple-rings} show that a non-\PI projectively simple algebra has a unique maximal ideal, namely the augmentation ideal.  A similar result holds for \ji algebras with nonzero Jacobson radical.

\begin{proposition}\label{thm:ji-not-semiprimitive-implies-finite-spec}
Let $A$ be a \ji $k$-algebra with nonzero Jacobson radical.  Then $A$ has only finitely many prime ideals.	
\end{proposition}

\noindent {\bf Proof.}  
We denote the Jacobson radical of the algebra $A$ by $\Jac(A)$.  Since $A$ is \ji and $\Jac(A)\neq (0)$, $A/\Jac(A)$ is finite dimensional.  It follows that there are only finitely many prime ideals in $A$ that properly contain $J(A)$ \cite[Theorem 2.3.9 (ii)]{Rowen:Ring-theory-Vol1}.  But every nonzero prime ideal in $A$ is necessarily maximal and thus primitive and so it contains $\Jac(A)$.  The result follows. 
\qed

The hypothesis that $\Jac(A) \neq (0)$ is necessary, as $\C[x]$ demonstrates.  Proposition \ref{thm:ji-not-semiprimitive-implies-finite-spec} motivates the following

\begin{question}\em{
If $A$ is a \fg non-\PI \ji $k$-algebra, does $A$ have only finitely many prime ideals? }	
\end{question}

It is known that a \PI \ji algebra is a \fg module over its center and the center is itself \ji \cite{FarinaPendergrass:A-Few-Properties-of-Just-Infinite}.  In stark contrast to the \PI case, a \ji algebra that doesn't satisfy a polynomial identity necessarily has a small center (assuming the ring isn't simple).  In fact, we next show that the center of any non-simple non-\PI \ji algebra is a finite dimensional field extension of the ground field.  This result is essentially a consequence of two lemmas of Farkas and Small \cite{FarkasSmall:Algebras-which-are-nearly}.  For a fixed positive integer $n$, we denote by $I_{n}(A)$ the two-sided ideal of the algebra $A$ generated by all specializations in $A$ of all polynomial identities of $n \times n$ matrices.

\begin{proposition}\label{thm:center-of-non-pi-ji}
Let $A$ be a non-\PI \ji $k$-algebra that is not simple.  Then $Z(A)$ is a finite dimensional field extension of $k$.
\end{proposition}

\noindent {\bf Proof.}
Since $A$ is \ji, it is prime \cite{FarinaPendergrass:A-Few-Properties-of-Just-Infinite}, and so $Z(A)$ is a domain.  Choose $0\neq z\in Z(A)$ and set $J = \bigcap_{d\in \N}Az^{d}$.  Since $A$ is \ji, the two-sided ideal $Az$ has finite codimension in $A$.  A result of Farkas and Small \cite[Lemma 1.1]{FarkasSmall:Algebras-which-are-nearly} gives $I_{n}(A) \subseteq J$ for some $n$.  Since $A$ doesn't satisfy a polynomial identity, $I_{n}(A)$ is a nonzero ideal of $A$, and so $J$ is nonzero as well.  It follows that $z$ is a unit of $A$ \cite[Lemma 1.2]{FarkasSmall:Algebras-which-are-nearly} and thus is invertible in $Z(A)$.  Thus $Z(A)$ is a field extension of $k$.  Next, since $A$ is not simple, there is a nonzero proper ideal $I$ of $A$.  Since $Z(A)$ is a field, $I\cap Z(A) = (0)$ and so the natural projection map $A \ra A/I$ induces an injection from $Z(A)$ to $A/I$, and the latter is finite dimensional.   
\qed

The fact that the center of a (non-simple) non-\PI \ji algebra is a finite dimensional field extension of $k$ shows that non-\PI \ji algebras are a reasonable generalization of central simple algebras.  In particular, showing that a \ji algebra is \sji can be regarded as an analog of the fact that a central simple algebra remains simple under extension of scalars.

\section{Central Closure}\label{sec:central-closure}

To prove our results we make use of the so-called \emph{extended center}.  This idea was largely inspired by some of the methods of Resco \cite{Resco:A-reduction-theorem-for-the-primitivity}, to which we owe a great debt.  

Let $A$ be a prime ring and consider the set of all right $A$-module homomorphisms $f\colon I \ra A$, where $I$ ranges over all nonzero two-sided ideals of $A$.  Martindale \cite{Martindale:Prime-rings-satisfying} shows how one can endow this set of maps with a useful algebra structure.  For the reader's convenience we include the full definition of Martindale's ring of quotients, but we refer to other sources \cite{:Handbook-of-algebra.-Vol., Lam:Lectures-on-modules-and-rings, Rowen:Ring-theory-Vol1} for the proof that this construction actually yields a ring with the stated properties.

\begin{definition}\label{def:martindale-ring-of-quotients}\em{
Let $A$ be a prime $k$-algebra.  The \emph{(right) Martindale ring of quotients of $A$}, denoted $\mart{A}$, consists of equivalence classes of pairs $(I, f)$ where $I \ideal A$, $I\neq (0)$, and $f\in \Hom_{A}(I_{A}, A_{A})$.  Here two pairs $(I,f)$, $(J, g)$ are defined to be equivalent if $f = g$ on the intersection $I\cap J$.  Addition and multiplication are given by}
\begin{align*}
(I, f) + (J, g) &= (I\cap J, f+g), \\
(I,f)\cdot(J,g) &= (JI, f\circ g).
\end{align*}
\end{definition}

\begin{definition}\em{
Let $A$ be a prime $k$-algebra.  The \emph{extended center of $A$}, written $\extcent{A}$, is defined to be $Z(\mart{A})$.  }
\end{definition}

\begin{remark}\em{
$\extcent{A}$ is also called the \emph{extended centroid} in the literature.  Hopefully our terminology won't cause confusion, as the centroid of a unital ring coincides with its center.  $\extcent{A}$ is a field extension of $k$ and $\extcent{A}\cap A = Z(A)$.  The \emph{central closure of $A$}, denoted $\centclose{A}$, is the $\extcent{A}$-linear subspace of $\mart{A}$ generated by $A$.  Lastly, $A$ is called \emph{centrally closed} if $\extcent{A} = Z(A)$, (equivalently, if $\centclose{A} = A$).  Note that the central closure of any ring is itself centrally closed.  }
\end{remark}

There is an entirely internal characterization of $\extcent{A}$ which bears mentioning.  The extended center of $A$ consists precisely of those pairs $(I, f)$ where $f\colon I \ra A$ is an $(A,A)$-bimodule homomorphism:
\[
\extcent{A} = \setmid{(I,f)}{(0) \neq I\ideal A, f\in \Hom_{A}(_{A}I_{A}, _{A}A_{A})}.
\]

Centrally closed prime algebras behave particularly well with respect to extension of scalars, and as a result we can show that a non-\PI \ji algebra is \sji if and only if it is centrally closed.

\begin{theorem}\label{thm:non-pi-sji-iff-centrally-closed}
Let $k$ be algebraically closed, and let $A$ be a non-\PI \ji algebra with $Z(A) = k$.  Then $A$ is \sji iff $A$ is centrally closed.	
\end{theorem}
\noindent {\bf Proof.}
When $A$ is simple (and hence centrally closed) the result follows from the well-known fact that central simple rings remain simple upon extension of scalars.  Thus we may assume that $A$ is not simple.  Suppose that $A$ is centrally closed.  Choose a nonzero ideal $I \ideal A\otimes_{k}K$ and set $J = I\cap A$, so $J \otimes_{k}K \subseteq I$.  Since $A$ is centrally closed, $J$ is a nonzero ideal of $A$ \cite[Lemma 3.4]{EricksonMartindale:Prime-nonassociative-algebras}.  Since 
\[
(A\otimes_{k}K)/(J\otimes_{k}K) \cong A/J \otimes_{k}K
\]
and the latter is finite dimensional over $K$, we see that $I$ has finite codimension in $A\otimes_{k}K$.

For the other direction, suppose that $A$ is \sji but not centrally closed.   Choose $(\varphi,I) \in \extcent{A} \setminus k$.  We may assume that $I$ is generated by a single element $a$, and so $\varphi$ is completely determined by $b:= \varphi(a)$.  Since $(\varphi,I) \notin k$, we may further assume that $b \neq \lambda a$ for any $\lambda \in k$.  Set $B := A\otimes_{k} k(t)$ where $t$ is an indeterminate.  We identify $B$ with the localization of $A[t]$ obtained by inverting all nonzero polynomials in $t$ with coefficients in $k$.  Consider the ideal $(a+tb)$ of $B$ generated by $a+tb$.  We will show that $A\cap (a+tb) = (0)$, so suppose there is some nonzero element $c \in A \cap (a+tb)$.  We may write $c = \sum_{i} r_{i}(a + tb) s_{i}$, with $r_{i}, s_{i} \in B$.  If $q(t)\in k[t]$ is a common right denominator for the $r_{i}$ and $s_{i}$, we see that    
\[
c q(t) = \sum_{i} f_{i}(t)(a+tb) g_{i}(t), \; \mbox{ for some } f_{i}, g_{i} \in A[t].   
\]
Now, writing $f_{i}(t) = \sum_{j} p_{ij} t^{j}$ and $g_{i}(t) = \sum_{l} q_{il} t^{l}$ where $p_{ij}, q_{il} \in A$, we see that
\begin{align*}
c q(t) &= \sum_{i,j,l} p_{ij} a q_{il} t^{j+l} +  p_{ij} b q_{il} t^{j+l +1} \\
 	&= \sum_{i,j,l} p_{ij} a q_{il} t^{j+l} + t \sum_{i,j,l} p_{ij} \varphi(a) q_{il} t^{j+l}.
\end{align*}
Set $u(t):= \sum_{i,j,l} p_{ij} a q_{il} t^{j+l} \in A[t]$, so
\[
c q(t) = u(t) + t \varphi(u(t)),
\]    	
where here $\varphi$ denotes the $(A,A)$-bimodule endomorphism of $A[t]$ defined by $\varphi(a) = b$ and $\varphi(t) = t$.  Note that for all $\lambda \in k$, 
\[
cq(\lambda) = u(\lambda) + \lambda \varphi(u(\lambda)) \in A,
\]
and so for all $\lambda \in k$ we have
\[
c (q(\lambda) - q(0)) = u(\lambda) + \lambda \varphi(u(\lambda)) - u(0).\tag{\dag}
\]

Let $V$ denote the $k$-span of $\setmid{u(\lambda)}{\lambda \in k}$, $W = V + kc$, and note that $W$ is a nonzero subspace of $A$ with $\dim_{k}(W) < \infty$.  Moreover, by equation $(\dag)$, 
\[
\varphi(W) \subseteq W,
\]
and so $\varphi$ is algebraic over $k$ by the Cayley-Hamilton theorem.  Since $k$ is algebraically closed, we see that $\varphi \in k$.  Thus $A \cap (a + tb) = (0)$ and hence $A$ embeds in $B/(a+tb)$.  However, $A$ is \sji so $B/(a+tb)$ is finite dimensional over $k(t)$, and thus $A$ is \PI.  This contradiction completes the proof.
\qed

The careful reader might note that we only used the non-\PI hypothesis in one direction of the preceding proof.  But note that if $A$ is \PI and centrally closed, then the center of $Q(A)$ is reduced to scalars.  Thus $Q(A)$ is finite dimensional over $k$ by Kaplansky's theorem and hence $A$ cannot be \ji.

If we remove the condition that $Z(A) = k$ in Theorem \ref{thm:non-pi-sji-iff-centrally-closed}, then it seems that to decide whether or not a non-\PI \ji algebra is \sji requires consideration of arbitrary field extensions $K/k$.  In fact, the next proposition shows that one only need consider the extended center.

\begin{proposition}\label{thm:sji-iff-tensor-with-extcent}
Let $A$ be a non-\PI \ji $k$-algebra, and let $\extcent{A}$ denote the extended center of $A$.  Then $A$ is \sji iff $A\otimes_{k}\extcent{A}$ is \ji over $\extcent{A}$.
\end{proposition}

\noindent {\bf Proof.}
One direction is trivial.  For the other, suppose that $A\otimes_{k}\extcent{A}$ is \ji over $\extcent{A}$.  Note that we have a surjective $k$-algebra homomorphism from $A\otimes_{k}\extcent{A}$ onto the central closure of $A$ (given by multiplication).  Since $A$ is \ji and not \PI, $\dim_{\extcent{A}}\centclose{A} = \infty$, and thus the above map is an isomorphism: $A\otimes_{k}\extcent{A} \cong \centclose{A}$.  We see that $\centclose{A}$ is \ji over $\extcent{A}$.
We claim then that $\centclose{A}$ is \sji over $\extcent{A}$.  If $A$ is simple then $\extcent{A} = Z(A)$, so $\centclose{A} = A$ is central simple over $\extcent{A}$ and thus $A$ is \sji.  If $A$ is not simple, then, as in the first part of the proof of Theorem \ref{thm:non-pi-sji-iff-centrally-closed}, it is easy to see that $\centclose{A}$ is \sji over $\extcent{A}$.  

Now, let $K/k$ be any field extension, and let $L$ denote a compositum (over $k$) of $K$ and $\extcent{A}$.  We then have that 
\[A\otimes_{k}L \cong (A\otimes_{k}\extcent{A})\otimes_{\extcent{A}}L
\] 
is \ji over $L$, and Lemma \ref{lem:sji-descends-to-subfields} shows that $A\otimes_{k}K$ is \ji over $K$, completing the proof.   
\qed

Theorem \ref{thm:non-pi-sji-iff-centrally-closed} and Proposition \ref{thm:sji-iff-tensor-with-extcent} give entirely internal characterizations of non-\PI \sji algebras.  The difficulty is that in practice it is often hard to compute the extended center.

\section{Noetherian \JI Algebras}\label{sec:noetherian-ji-rings}

In this section we prove that a non-\PI right noetherian \ji algebra over an algebraically closed field is \sji.
To obtain this result, we use central closure.

\begin{proposition}
Let $A$ be a right noetherian non-\PI \ji algebra with center $k$.  Then $\extcent{A}$ is an algebraic extension of $k$.
\end{proposition}

\noindent {\bf Proof.}
Since $A$ is right Goldie, the right maximal ring of quotients of $A$ coincides with its right Goldie quotient ring, which we denote by $Q(A)$.  Since $A$ is prime, $\extcent{A} = Z(Q(A))$ \cite[Prop 14.17]{Lam:A-first-course-in-noncommutative}.  Choose a nonzero element $t \in \extcent{A}$ and let $B$ be the subalgebra of $Q(A)$ generated by $A$ and $t$.  Let $p = p(t)$ be any nonzero polynomial in $t$ with coefficients in $k$ and let $\pi : B \ra B/pB$ denote the natural projection map.  There are nonzero elements $a,b\in A$ with $b$ regular such that $p = ab^{-1}$ and so $pb = a$ in $B$.  Hence $pB \cap A\neq (0)$.  Moreover, 
\[
B/pB \cong \left(\frac{A}{pB \cap A}\right) \set{\overline{t}}, \; \mbox{ where } \overline{t} = t + pB.
\]
Thus every element of $B/pB$ can be written as a polynomial in $\overline{t}$ of degree $< \deg(p)$ with coefficients in $A/(pB \cap A)$.  Since $A$ is \ji, this shows that $\dim_{k}(B/pB) < \infty$.  Set $n = \dim_{k}(B/pB)$.  By a result of Farkas and Small \cite[Lemma 1.1]{FarkasSmall:Algebras-which-are-nearly}, $I_n(B) \subseteq \tcap_{d\geq 1}p^{d}B$, and since $A$, and hence also $B$, is not \PI, $I_{n}(B) \neq (0)$.

Choose any nonzero regular element $$y \in \tcap_{d\geq 1} p^{d}B$$ (such a regular element exists because $B$ is prime right Goldie).  Since $p$ is central, we can then find regular elements $a_{1}, a_{2}, \ldots \in B$ with 
\[
y = a_{d}p^{d}, \;\mbox{ for all } d\in \N,
\]
and since $p$ is regular in $B$, we have $a_{d} = a_{d+1}p$ for all $d$.  In particular, $a_{d}B \subseteq a_{d+1}B$.  Since $B$ is right noetherian, the chain of right ideals 
\[
a_{1}B \subseteq a_{2}B \subseteq \cdots
\]
must terminate, so $a_{i}B = a_{i+1}B$ for some $i$.  In particular, $a_{i}B = a_{i+1}B = a_{i+1}p B$, and since $a_{i+1}$ is regular, the second equality shows that $p$ has a right inverse in $B$.  Since $B$ is right noetherian, $p$ is a unit in $B$.
  
Let $\mf{m}\in \mspec(B)$ be any nonzero maximal ideal of $B$.  We have a $k$-algebra isomorphism
\[
\frac{B}{ \mf{m} } \cong \left(\frac{A}{\mf{m}\cap A}\right)\set{\overline{t}}, \mbox{ where } \overline{t} = t + \mf{m}.
\] 
Since $A$ is \ji, $A/(\mf{m}\cap A)$ is finite dimensional over $k$, and since $\overline{t}$ is central, $B/\mf{m}$ is a simple \PI ring, and thus a finite module over its center, $Z$, by Kaplansky's theorem.  Since $B/\mf{m}$ is clearly an affine $k$-algebra, the Artin-Tate lemma then shows that $Z$ is affine as well.  But $Z$ is a field, and thus $\dim_{k}(Z) < \infty$.  This shows that $\dim_{k}(B/\mf{m}) < \infty$.

Finally, since $B/\mf{m} \cong (A/(\mf{m}\cap A))\set{\overline{t}}$ is finite dimensional, we see that $\overline{t}$ is algebraic over $k$, and we claim that this implies that $t$ is algebraic over $k$ as well.  To see this, suppose that $t$ were transcendental.  Since any nonzero polynomial $p(t)\in k\set{t}$ is a unit in $B$, $B$ contains the field $k(t)$, which has transcendence degree $1$.  Since $\mf{m}\neq B$, $\mf{m} \cap k(t) = (0)$, and so $B/\mf{m}$ contains an isomorphic copy of $k(t)$, namely $k(\overline{t})$. Since $\overline{t}$ is algebraic over $k$, this is a contradiction.
\qed

Since being centrally closed and being \sji are equivalent properties for non-\PI \ji algebras, we obtain as a corollary Theorem \ref{thm:non-pi-noetherian-ji-implies-sji}.

\begin{corollary}\label{cor:noetherian-non-pi-ji-implies-sji}
Let $k$ be an algebraically closed field, and let $A$ be a right noetherian non-\PI \ji algebra with $Z(A) = k$.  Then $A$ is \sji.
\end{corollary}

We can obtain further stability results by combining the above ideas with the Nullstellensatz.

\begin{proposition}\label{thm:primitive-null-implies-sji}
Let $k$ be an algebraically closed field and let $A$ be a primitive non-\PI \ji $k$-algebra.  If $A$ satisfies the Nullstellensatz, then $A$ is \sji.
\end{proposition}
 
\noindent {\bf Proof.}
By a result of Martindale \cite{Martindale:Lie-isomorphisms-of-prime-rings}, the extended center $C(A)$ embeds in $\End_{A}(M)$ for any faithful simple $A$-module $M$.  Since $A$ satisfies the Nullstellensatz, this implies that $C(A) = k$, so $A$ is centrally closed and the result follows from Theorem \ref{thm:non-pi-sji-iff-centrally-closed}.
\qed

\section{\PI \JI Algebras}\label{sec:pi-ji-rings}

In this section, we characterize \PI \sji rings over algebraically closed fields.  We begin by considering the commutative case.

\begin{lemma}\label{lem:noetherian-domains-of-trdeg-1-are-ji}
Let $k$ be an algebraically closed field and let $A$ be a commutative noetherian domain such that $Q(A)$ is a \fg field extension of $k$ with $\trdeg(Q(A)) \leq 1$, then $A$ is \ji.	
\end{lemma}
\noindent {\bf Proof.}
We denote the (classical) Krull dimension of a ring $A$ by $\Krull(A)$.   Since $A$ is a domain with $\Krull(A) \leq 1$, $\Krull(A/P) = 0$ for every nonzero $P\in\spec(A)$.  Thus $A/P$ is a field extension of $k$.  If $A$ is not \ji, then, since $A$ is noetherian, we may choose an ideal $I$ maximal with respect to $\dim_{k}(A/I) = \infty$.  It is easy to see that $I$ is then prime and it suffices to show that $A/P = k$ for every nonzero prime ideal $P$ of $A$.  Fix a nonzero $P\in \spec(A)$ and suppose there is some element $x + P\in A/P$ which is transcendental over $k$.  Since $k$ is algebraically closed and $A$ is a domain, every nonzero element of $P$ is transcendental over $k$.  Fix a nonzero $t\in P$.  Then $x$ and $t$ are algebraically independent over $k$ since $x$ is regular mod $P$, and so $\trdeg(Q(A)) \geq 2$, a contradiction.  Thus we see that $A/P$ is an algebraic extension of $k$.  Since $k$ is algebraically closed field, we see that $A/P\cong k$.
\qed

To characterize commutative \sji rings, we first require another lemma.

\begin{lemma}\label{lem:sji-descends-to-subfields}
Let $A$ be a $k$-algebra and let $K/k$ be a field extension.  If $A\otimes_{k}K$ is \ji over $K$, then $A$ is \ji over $k$.
\end{lemma}

\noindent {\bf Proof.}
Choose a nonzero ideal $I \ideal A$.  Then $I\otimes_{k}K$ is a nonzero two-sided ideal of $A\otimes_{k}K$ with $(A\otimes_{k}K)/(I\otimes_{k}K) \cong (A/I)\otimes_{k}K$.  Taking dimensions then yields
\begin{align*}
\dim_{k}(A/I) &= \dim_{K}( (A/I)\otimes_{k} K ) \\
	&= \dim_{K} ( (A\otimes_{k}K)/(I\otimes_{k}K) ) < \infty. \qed
\end{align*}

\begin{proposition}\label{thm:commutative-sji-characterization}
Let $k$ be an algebraically closed field and let $A$ be a commutative \ji $k$-algebra.  Then $A$ is \sji iff $Q(A)$ is a \fg field extension of $k$ of transcendence degree at most $1$.
\end{proposition}

\noindent {\bf Proof.}  Suppose first that $A$ is \sji.  
Set $F = Q(A)$, the quotient field of $A$.  If $A$ is \sji, then $A\otimes_k K$ is noetherian for every field extension $K/k$ since a commutative just infinite ring is necessarily commutative.  Since $F\otimes_{k} K$ is a localization of $A\otimes_k K$, it is also noetherian for every field extension $K/k$.  By a result of V\'{a}mos \cite{Vamos:On-the-minimal-prime-ideal}, $F$ is a \fg field extension of $k$, and it suffices to show that $\trdeg_{k}(F) \leq 1$.  Note that $\Krull(A\otimes_{k} F) \leq 1$ since $A\otimes_{k} F$ is a \ji $F$-algebra.  Since $F\otimes_{k} F$ is a localization of $A\otimes_{k} F$, $F\otimes_{k} F$ is a \ji $F$-algebra as well, and so $\Krull(F\otimes_{k} F) \leq 1$, and thus $\trdeg_{k}(F)\leq 1$ \cite[6.4.8 \& 6.6.17]{McConnellRobson:Noncommutative-Noetherian-rings}.

Next suppose that $F$ is a \fg field extension of $k$ with $\trdeg_{k}(F) \leq 1$.  We claim that $A$ is \sji.  Let $K/k$ be a field extension, and let $\overline{K}$ denote the algebraic closure of $K$.  Since $F$ is \fg, there is a \fg purely transcendental extension $E/k$ with $E\subseteq F$ and $[F:E]<\infty$. Then $\trdeg_k(E)\le 1$ and so the field of fractions of 
$E\otimes_k \overline{K}$ is a purely transcendental extension of $\overline{K}$ of transcendence degree at most one.  But $Q(F\otimes_{k} \overline{K})$ is a finite extension of $E\otimes_k \overline{K}$ and thus also has transcendence degree at most one.  Since $Q(A\otimes_k \overline{K})=Q(F\otimes_k \overline{K})$,  $A\otimes_{k} \overline{K}$ has Krull dimension at most $1$ over $\overline{K}$.  Applying Lemma \ref{lem:noetherian-domains-of-trdeg-1-are-ji}, we see that $A\otimes_{k} \overline{K}$ is \ji over $\overline{K}$.  Lemma \ref{lem:sji-descends-to-subfields} then implies that $A\otimes_{k} K$ is \ji over $K$, completing the proof.   
\qed

\begin{corollary}
Affine commutative \ji algebras over algebraically closed fields are \sji.	
\end{corollary}

Suppose now that $A$ is a \PI \ji algebra over an algebraically closed field $k$ such that $Q(Z(A))$ is a \fg field extension of $k$ with $\trdeg_{k}(Q(Z(A))) \leq 1$.  Then $A$ is a \fg module over $Z(A)$, which is itself \ji \cite[Corollary 2]{FarinaPendergrass:A-Few-Properties-of-Just-Infinite}.  By Proposition \ref{thm:commutative-sji-characterization}, $Z(A)$ is \sji.

If $K/k$ is a field extension, then $Z(A)\otimes_{k}K$ is a (noetherian) \ji $K$-algebra.  Since $A\otimes_{k}K$ is a \fg $Z(A)\otimes_{k}K$-module, $A\otimes_{k}K$ is also noetherian.  Choose a nonzero two-sided ideal $I\ideal A\otimes_{k}K$.  As $A\otimes_{k}K$ is prime noetherian, $I$ contains a nonzero regular element $x$.  Since $A\otimes_{k}K$ is a finite $Z(A)\otimes_{k}K$-module, $x$ is algebraic over $Z(A)\otimes_{k}K$.  If 
\[
x^{t} + z_{t-1}x^{t-1} + \ldots + z_{0} = 0, \quad z_{j}\in Z(A)\otimes_{k}K
\]      
has $t$ minimal, then $z_{0} \neq 0$.  It follows that $z_{0}\in I\cap (Z(A)\otimes_{k}K)$, and so $I \cap (Z(A)\otimes_{k}K) \neq (0)$.  Finally,
$(A\otimes_{k}K)/I$ is a \fg module over $(Z\otimes_{k}K)/(I\cap (Z\otimes_{k}K))$, hence $A\otimes_{k}K$ is \ji over $K$.  To summarize, we have the following

\begin{theorem}\label{thm:pi-sji-characterization}
Let $k$ be an algebraically closed field, and let $A$ be a \PI \ji $k$-algebra.  Then $A$ is \sji iff $Q(Z(A))$ is a \fg field extension of $k$ with $\trdeg_{k}(Q(Z(A))) \leq 1$.
\end{theorem}

\begin{corollary}\label{cor:affine-pi-over-alg-closed-is-sji}
Affine \PI \ji algebras over algebraically closed fields are \sji.	
\end{corollary}

\section{Algebras over Uncountable Fields}\label{sec:uncountable-fields}

In this section we show that countably generated \ji algebras over uncountable algebraically closed fields are \sji, proving Theorem \ref{thm:cg-ji-over-big-fields-are-sji}.  We first consider the \PI case.  

\begin{proposition}\label{thm:cg-comm-ji-rings-over-uncountable-are-sji}
Let $k$ be an uncountable algebraically closed field and let $A$ be a countably generated \PI \ji $k$-algebra.  Then $A$ is \sji.	
\end{proposition}
\noindent {\bf Proof.}
We begin by showing we can reduce to the commutative case.  Since $A$ is \ji and \PI, we have that $A$ is a finite module over its center and that its center is \ji \cite{FarinaPendergrass:A-Few-Properties-of-Just-Infinite}.  Moreover, since an algebra is countably generated iff it is countable dimensional, we see that $Z(A)$ is countably generated as well.  By Theorem \ref{thm:pi-sji-characterization}, we may then assume that $A$ is commutative.  

Note that $\Krull(A) \leq 1$ since $A$ is \ji.  In fact, $\Krull(A) = 1$ because there are no proper field extensions of $k$ which are countably generated as algebras.  By Proposition \ref{thm:commutative-sji-characterization} it is enough to show that $Q(A)$ is a \fg field extension of $k$ of transcendence degree at most $1$.  $A$ is noetherian, so $A\otimes_{k} Q(A)$ is noetherian \cite[Theorem 1.2]{Bell:Noetherian-algebras-over}.  Thus $Q(A)\otimes_{k} Q(A)$ is noetherian since it is a localization of $A\otimes_{k}Q(A)$.  It follows from a result of V\'{a}mos \cite{Vamos:On-the-minimal-prime-ideal} that $Q(A)$ is a \fg field extension of $k$.  

Let $R$ be any affine commutative $k$-algebra with $Q(R) = Q(A)$.  Since $A$ is countably generated, we may choose elements $a_{1}b_{1}^{-1}, a_{2}b_{2}^{-1}, \ldots$, with $a_{i}, b_{i} \in R$, which generate $A$ as a $k$-algebra.  Set $B = RS^{-1}$, where $S$ is the multiplicatively closed subset of $R$ generated by the set $\setmid{b_{i}}{i \in \N}$.  Note that $B$ satisfies the Nullstellensatz since it is countably generated and $k$ is uncountable.  Since $k$ is algebraically closed, $B/M \cong k$ for some maximal ideal $M$ of $B$.  Set $P = M \cap R$.  Then $P$ is a maximal ideal of $R$.  Note also that $R_{P}$, the localization of $R$ at $P$, contains $B$ and thus also contains $A$.  Hence $\Krull(A) = \trdeg_{k}(A)$ \cite[Theorem 1.5]{GilmerNashierNichols:The-prime-spectra-of-subalgebras}.
\qed

Next we dispose of the non-\PI case.  The techniques are rather different, and we need a preliminary result which may be of independent interest.

\begin{proposition}\label{thm:extended-center-is-countable-dimensional}
Let $A$ be a countably generated non-\PI \ji $k$-algebra.  Then $\extcent{A}$ is a countable dimensional field extension of $k$.
\end{proposition}

\noindent {\bf Proof.}
If $A$ is simple, then $Z(A)$ is a field.  But $A$ is countable dimensional over $k$, so $Z(A)$ is countable dimensional over $k$ as well, and in a simple ring, $\extcent{A} = Z(A)$.

Now suppose that $A$ is not simple.  Since $A$ doesn't satisfy any polynomial identity, each $I_{n}(A)$ is nonzero.  Moreover, $A$ is \ji, and so every nonzero ideal $I \ideal A$ contains some $I_{n}(A)$.  Suppose there is an uncountable collection of elements 
\[
\setmid{ (I_{\alpha}, f_{\alpha}) }{ \alpha \in \mathcal{I} } \subseteq \extcent{A}
\]  
which is $k$-linearly independent.  Note that if $(I,f) \in \extcent{A}$, then 
\[
(I,f) = (I_{n}(A), f|_{I_{n}(A)}),
\] 
so every bimodule map is determined by it's action on some $I_{n}$.  Thus for each $\alpha$, we may assume that $(I_{\alpha}, f_{\alpha}) = (I_n(A), f_{\alpha})$ for some $n$.  Since there are only countably many $I_{n}(A)$, there is some $d\in \N$ for which the uncountable set $\setmid{ (I_{d}(A), f_{\alpha}) }{ \alpha\in \mathcal{I} }$ is $k$-linearly independent.  Next, choose any nonzero element $a\in I_{d}(A)$ and note that by replacing $I_{d}(A)$ with $AaA$, we may assume that there is an uncountable $k$-linearly independent set of the form $\setmid{ (AaA, f_{\alpha}) }{\alpha \in \mathcal{I}}$.  Being an $(A,A)$-bimodule map, each $f_{\alpha}$ is completely determined by $f_{\alpha}(a)$.  However, since $A$ is a countably generated $k$-algebra, $\dim_{k}(A)$ is countable, and it follows that the set 
\[
\setmid{ f_{\alpha}(a) }{ \alpha \in \mathcal{I}}
\] 
is $k$-linearly dependent.  This proves the theorem.  
\qed

\begin{remark}\em{
Though Proposition \ref{thm:extended-center-is-countable-dimensional} requires $A$ to be non-\PI and \ji, the conclusion still holds for any countably generated prime algebra with the property that there is a countable set of nonzero elements such that every nonzero ideal contains at least one of these elements; this is often referred to as having a \emph{countable separating set} in the literature.  Note that, for example, any primitive homomorphic image of an enveloping algebra of a finite dimensional Lie algebra has a countable separating set, as does any countably generated prime non-\PI right Goldie ring of $\GK$ $2$.  We also remark that the preceding proof, with the obvious modifications, yields a generalization of \cite[Theorem 1]{RowenSmall:Primitive-ideals-of-algebras}.	}
\end{remark}

\begin{corollary}\label{cor:cg-non-PI-ji-over-big-fields-are-sji}
Let $k$ be an uncountable algebraically closed field and let $A$ be a countably generated non-\PI \ji $k$-algebra.  Then $A$ is \sji.	
\end{corollary}
\noindent {\bf Proof.}
If $A$ is simple then $\extcent{A} = k$ and the result follows from Theorem \ref{thm:non-pi-sji-iff-centrally-closed}.  Otherwise $A$ is not simple, and so $\extcent{A}$ is a countable dimensional field extension of $k$ by Proposition \ref{thm:extended-center-is-countable-dimensional}.  Since $k$ is uncountable, $\extcent{A}$ must be algebraic over $k$.  But $k$ is algebraically closed, so $A$ is centrally closed, and the result follows from Theorem \ref{thm:non-pi-sji-iff-centrally-closed}.
\qed

Combining Proposition \ref{thm:cg-comm-ji-rings-over-uncountable-are-sji} and Corollary \ref{cor:cg-non-PI-ji-over-big-fields-are-sji}, we obtain Theorem \ref{thm:cg-ji-over-big-fields-are-sji}.

\section{An Example}\label{sec:examples}

In this short section we prove Theorem \ref{thm:affine-ji-non-sji-rings-over-countable-fields-exist}. 
We construct an interesting example showing that the right noetherian hypothesis in Theorem \ref{thm:non-pi-noetherian-ji-implies-sji} and the uncountable hypothesis in Theorem \ref{thm:cg-ji-over-big-fields-are-sji} are necessary.

Let $k$ be a countable algebraically closed field.  There exists a simple Ore domain $R$, affine over $k$, with $\trdeg_{k}(Z(R)) = \infty$ \cite{Irving:Finitely-generated-simple}.  Thus $R$ is affine, simple, and infinite dimensional, and so $R$ cannot satisfy a polynomial identity.  The ring $R$ is obtained as a skew laurent extension of an Ore domain with respect to an automorphism of infinite order \cite{Irving:Finitely-generated-simple}.  In particular, $R$ is not a division ring.  Fix a nonzero element $x\in R$ which is not a unit, and set $A = k + xR$.  $A$ is affine \cite[Proposition 2]{Resco:Affine-domains-of-finite}, and since $xR$ is simple as a (non-unital) ring, $A$ is \ji.    As $R$ is not \PI, neither is $A$, and so we see that $Z(A)$ is a finite dimensional field extension of $k$ by Proposition \ref{thm:center-of-non-pi-ji}.  But note that $Z(A) \cap xR = (0)$ because $xR$ is a proper right ideal of $R$.  Hence $Z(A) = k$.  Finally, it is easy to see that $A$ is right Ore, and so the extended center of $A$ coincides with the center of it's right Goldie quotient ring.  Since $Q(A) = R$, we see that $\extcent{A} = Z(R)$.  Thus $A$ is not centrally closed, and so $A$ is not \sji by Theorem \ref{thm:non-pi-sji-iff-centrally-closed}.

\section*{Acknowledgments}

We thank Lance Small for many helpful conversations.

\bibliographystyle{abbrv}

\end{document}